\documentclass[a4paper,12pt]{opelsarticle}
\usepackage{amsmath,amsfonts,graphicx,tikz,pgfplots, tikz-3dplot, pdfsync}
\usepackage[colorlinks=true, citecolor=red]{hyperref}
\usepackage{fullpage}
\usetikzlibrary{patterns}
\usepackage{bm} 

\def\eq#1{\Blue{\begin{eqnarray}#1\end{eqnarray}}}

\DeclareGraphicsExtensions{.pdf, .jpg, .tif, .png}

\def\e{{\rm e}}

\def\R{{\mathbb R}}

\def\E{{\bf E}}
\def\R{{{\bf R}}}
\def\p{\partial}
\def\n{\nabla}

\def\d{\mathrm d}

\newtheorem{remark}{Remark}

\def\n{\nabla}
\def\E{{\mathbb E}}
\def\p{\partial}
\def\ds{\displaystyle }

\def\RR{{\mathbf{R}}}

\def\R{{\mathbb R}}

\def\U{{\mathcal U}}

\def\vr{{\mathbf r}}
\def\vu{{\mathbf u}}
\def\vv{{\mathbf v}}
\def\vv{{\mathbf v}}

\def\vn{{\mathbf n}}

\def\vI{{\mathbf I}}
\def\vB{{\mathbf B}}

\def\vN{{\mathbf N}}

\def\vZ{{\mathbf Z}}
\def\vW{{\mathbf W}}
\def\vQ{{\mathbf Q}}
\def\vX{{\mathbf X}}

\def\vx{{\mathbf x}}

\def\vb{{\mathbf b}}

\def\vA{{\mathbf A}}

\def\p{\partial}
\def\n{\nabla}

\def\d{\mathrm d}

\def\eq#1{{\begin{eqnarray}#1\end{eqnarray}}}

\begin{document}

\title{Optimal Quota for a Multi-species Fishing Models}
\author{
Olivier Pironneau\footnote{\emph{olivier.pironneau@sorbonne-universite.fr }, LJLL, Sorbonne Universit\'e, Paris, France.} \\~\\ 
{\it Dedicated to E. Polak for his 95$^{th}$ birthday}
}

\begin{frontmatter}
\begin{abstract}
A Stochastic Control Problem can be solved by Dynamic Programming or Distributed Optimal Control with the Kolmogorov equation for the probability density of the Markov process of the problem. It can be solved also with Supervised Learning.
We shall compare these two classes of methods for the control of fisheries.

Fishing quotas are unpleasant but efficient to control the productivity of a fishing site. A popular model has a vector-valued stochastic differential equation for the biomass of the different species. Optimization of quota will be obtained by a gradient method applied to the least square difference with an ideal state weighted by the probability density of the biomasses. Alternatively a deep neural network which preserves the  Markov property of the problem can be trained with a stochastic gradient algorithm.  The model is extended to distributed fishing sites and biomass is stabilized by adjusting the quota to its time derivative.
 \end{abstract}
\end{frontmatter}


\setcounter{tocdepth}{1}



\section*{Introduction}
The increasing need for food has led to over fishing everywhere.
To avoid extinction one must measure or  model the biomass and experiment with various ways to control it.
The mathematics of population dynamics are old (see Verhulst \cite{verhulst}). For  competing species (fish included)  Volterra \cite{lotka} introduced the  logistic predator-prey model in 1931. Since then, the model has been extended  and used by many (see for instance \cite{ALL},\cite{haddon} and \cite{mcglade}) and the literature is enormous.  
For fisheries Schaefer \cite{schaefer} introduced an effort function $E(t)$ -- conveniently representing the number of fishing boats at sea-- and a catchability coefficient $q$ for each class of boats.  In  \cite{MP18} an extension relating the fishing effort to the market price $p$ of fish is analyzed.  

Multi-species models are  straightforward vector generalizations of single species models, however their mathematical analysis and computer solutions are much harder.  The special case of a single species with different aging groups is usually analyzed by standard population dynamics arguments (see ``aged structured models'' in \cite{haddon}). Nevertheless, the complexity of the modeling can be grasped from \cite{king}, p73. 

 The Mathematical  literature on fishing quota is scarce \cite{Punt}. In \cite{katsukawa},\cite{WEI},\cite{DAN} the models are either too simple or analyzed in general terms for profitability and preservation  without numerical simulations. 
 
Our purpose, in this article, is to show what stochastic optimization can offer to fisheries. We heave no competence to discuss the accuracy of the models in practice.
 
 In \cite{PAOP} Supervised Learning was shown to be  efficient to calibrate the parameters of the fishing model of \cite{MP18}. In \cite{MLGPOP2} a stochastic control problem was derived for the computation of optimal quotas, a solution by Supervised Learning was proposed and compared to standard  stochastic control solutions using the Hamilton-Jacobi-Bellman equations (HJB).  

In this article  we compare a Distributed Control Method (an alternative to HJB)  to a new deep neural network which is an interesting modification (due to P. Bras \cite{PBGP}) of  the one used in \cite{MLGPOP2}.  A final remark about ``common sense control'' is made.

In the last section the model of \cite{pagespironneau} and \cite{MLGPOP2} is  extended to distributed fishing sites and solved numerically by ``common sense control''  for the Atlantic ocean facing Senegal. Some references to multi-sites models are available in \cite{moussaoui} and for open sea models in\cite{katsukawa}.

\section{The Single Fishing Site Model}

In simple situations, depleting of a sight due to fishing is proportional to the fish biomass $B$ and related to the fishing effort $E$ (the number of boats at sea) by
 \begin{equation}\label{CE1}
    \frac{\d B}{\d t}= B(r  - \kappa B) - q B E,
\end{equation}
Here $r$ is the natural birth minus death rate, $r/\kappa$ is the capacity of the site and $q$ is the catchability.
The rate of the fishing effort is proportional to the difference between profit $p B$ -- where $p$ is the price of fish -- and the cost $c$ of operating a fishing boat:
 \begin{equation}\label{CE2}
  \frac{dE}{dt}= p q B E  - c E .
\end{equation}
When the market is liquid the price adjusts daily to balance supply $qB E$ and demand $D(p)$, taken here inversely proportional to $1 + b p$ with $b$ fitted from past data. Thence a value for $p$ is found and the model can be rescaled to
\eq{\label{CE3}
\frac{\d B}{\d t}= B (r  - \kappa B - q E) ,~~
\frac{\d E}{\d t}=  a-(q B +c)E,~B(0)=B_0,~E(0)=E_0.
}
The model is easily extended to multi-species including a fishing quota $Q_i<q$ on each species $i=1,...,d$ and noise:
\eq{\label{CE4}&
\d\vB_t = \vB\star\left[\left(\vr  - \underline{\bm\kappa} \vB   -  \vQ  E\right)\d t + \underline{\bm\sigma}\d\vW_t\right] ,~~ &\vB(0)=\vB^0+\underline{\bm\sigma}"\vN_0^1,
\cr&
\d E_t=  \left(a-(\vB: \vQ  +c) E\right)\d t +E\underline{\bm\sigma}'\d\vW'_t, ~~&E(0)=E^0+\sigma N_{0,1}.
}
where $\underline{\bm\kappa}$ is the capacity matrix, $A\star B$ is the vector of component $A_i B_i$ and where  $A:B$ is the sum of $A_iB_i$.
$\vW$, $\vW'$, $\vN_0^1$,$N_{0,1}$ are Gaussian noises and $\underline{\bm\sigma}$, $\underline{\bm\sigma}'$, $\underline{\bm\sigma}''$, $\sigma$ are the variance-correlation matrices and variance coefficient. Note that the sign of $\bm\kappa_{ij}$ indicates whether species $i$ eat or is eaten by species $j$.    Noises are mathematical representations of the uncertainties on the parameters and on the model.  
\section{Identification of Coefficients}

Consider for simplicity a single species in absence of noise and assuming that $q$ is known; then $z:=[r,\kappa,a,c]$ must be identified.  The easiest is to choose two dates $t_1,t_2$ and measure $Z^d:=[X(t_1), E(t_1), X(t_2), E(t_2)]$.  It amounts to counting the number of boats at sea and how much fish were caught, on two different days.
Surprisingly, a root finding algorithm like \texttt{broyden1}  (from the Python library \texttt{scipy}) works very well \cite{PAOP} on synthetic data (i.e. choose a set $z_0$ to compute $Z(z_0)$, then invert numerically the mapping $z\mapsto Z$). The same can be achieved by least squares on the gap between the current state $Z$ and an ideal state $Z^d$. With noise, $\E$ being the expected value, one must solve.
\[
\min_z \E[|Z-Z^d|^2] ~:~\hbox{ subject to \eqref{CE3}}.
\]
This is a hard nonlinear stochastic optimization problem which is most likely not well posed before discretization because it uses discrete times.  Using Dynamic Programming  and Ito's formula to establish the optimality conditions a numerical solution requires to solve at each iteration of the optimization algorithm two partial differential equations \cite{MLGPOP2}; so it also expensive.

An easier solution can be obtained with a neural network to represent $Z(z)\mapsto z$ (two inner layers of 50 neurons + \texttt{ReLU} seem appropriate) and train the network as follows:
\begin{enumerate}
\item Prepare M synthetic solutions $\{Z(z^j)\}_1^M$ by solving \eqref{CE3}.
\item Train the network with the samples inputs $\{Z(z^j)\}_1^M$ and outputs $\{z^j\}_1^M$, using a least-square loss.
\end{enumerate}
Table \ref{tab2} shows typical results  for 3 values of the noise ($\sigma''=\sigma'=\sigma$) computed with a Neural Network made of 2 hidden layers with 100 neurons each and compared with  Dynamic Programming solutions.
{\small
\begin{table}[htp]
\caption{Identification of the parameters with a Neural Network in the random case.  Relative errors from the noiseless solution of \eqref{CE3} with $B^0=0.1$, $E^0=0.1$,  $r=2, \kappa=1$, $c=1$, $a=1.1$ using 1000  solutions (samples) of \eqref{CE3} at $t_1=1/14$, $t_2=1$, $q=1$  and 200 iterations (epochs). Dynamic Programming minimized the criteria around $5\cdot 10^{-3}$ with gradient values around $10^{-6}$.}
\begin{center}
\begin{tabular}{|c||c|c|c|c||c|c|c|c|}
$\sigma$&$r_{NN}$ &$\kappa_{NN}$ &$c_{NN}$ &$a_{NN}$& $r_{DP}$ &$\kappa_{DP}$ &$c_{DP}$ &$a_{DP}$\cr
0.01 & 1.99$\pm$ 0.09& 1.01$\pm$ 0.30& 0.97$\pm$ 0.06 &1.09$\pm$ 0.04&
1.95 & 0.74 &1.47 & 1.46\cr
0.125& 2.04$\pm$ 0.11&  1.13$\pm$ 0.20& 1.14$\pm$ 0.16 &1.29$\pm$ 0.10&
1.76&1.027&0.65&0.85\cr
0.25& 1.97$\pm$ 0.16& 1.03$\pm$ 0.34& 0.90$\pm$ 0.23& 1.15$\pm$ 0.15&
1.80 &  1.5 & 1.5 & 1.37
\cr\end{tabular}
\end{center}
\label{tab2}
\end{table}}

Supervise Learning gives a better solution in this case.

\section{Fishing Quotas}
Consider the problem of finding a suitable quota $\vQ(t)$ given to each fisherman for each species.   
We assume that $\vQ_i<q$ for all $i$, otherwise the fishermen are not affected and the quota is theoretical. Accordingly the total daily catch will be less than $\vB\star\vQ E$; this then is a global quota. Let $\vu:=\vQ E$;  searching for $\vu$ instead of $\vQ$ no longer requires the knowledge of $E$ and   $\vu\le \vu_M$ means that a global quota of $\vB\star\vu_M$ is imposed. To translate it at the fisherman level requires an estimate of $E$ (the number of boats at sea) before declaring the quota. As illegal fishing is hard to estimate, randomness in the model is welcome!

 Mathematically we may solve
\eq{\label{reform}&
\min_{\vu\in \U}\Big\{& \bar J:=\int_0^T\E\left[|\vB(t)-\vB^d|^2\d t  - \bm\alpha \cdot \vu + \bm\beta\cdot[\vu]_t^{0,T} \right] \d t ~: 
\cr&&
\d\vB_t=\vB \star\left[(\vr-\vu  - \underline{\bm\kappa}\vB)\d t +\underline{\bm\sigma}\d\vW_t \right],~~\vB(0)=\vB^0 + \underline{\bm\sigma}'\vN_0^1\Big\}.
}
The expectation is with respect to the laws on $\vW_t$ and $\vB^0$.
To preserve the Markovian feature of the problem we assume that $\vu$ is a deterministic function $\vx$ and $t$.
Also $\U=\{\vu~:~\vu_j\in[u_m,u_M],~j=1..d\}$. The quadratic variation is,
 \[
	 [\vu]^{0,T}_t = \lim_{\|P\|} \sum_{k=1..}^{t_k<t}|\vu_{t_k}-\vu_{t_{k-1}}|^2.
\]
where $P$ ranges over partitions of the interval $(0,t)=\cup_k(t_{k-1},t_k)\subset(0,t),~t<T$ and the limit is in probability when $\max_k|t_k-t_{k-1}|\to 0$.
Here It\^o calculus \cite{BIC} tells us that:
 \[
 \E[\vu]^{0,T}_t = \int_0^t \E[|\bm{\underline\sigma}\vB_t\cdot\n_\vB \vu|^2] \d t.
 \]
The term $\bm\alpha\cdot\vu$ encourages large quotas and represents the political cost of constraining the fishermen with small quotas; the term with $\bm\beta$ is added to prevent large oscillations of $\vu$ from one day to the next.
In \cite{MLGPOP2} it is shown that the problem is well posed. A solution exists but it may not be unique. Three numerical methods for solutions have been analyzed in \cite{MLGPOP2}: Stochastic Dynamic Programming, Hamilton-Jacobi-Bellman dynamic programming (HJB), and using Deep Neural Networks (DNN).  Here we present a modified DNN proposed in \cite{PBGP} and compare the results with the solution of the (equivalent) distributed control problem using Kolmogorov's forward equation for the probability density of $\vB$.

\subsection{The Distributed Control Problem}
Assume for clarity that $\underline{\bm\sigma}=\sigma\vI,\sigma$ constant and $\bm\beta_i=\beta$ for all $i$. The Kolmogorov equation for  $\rho(\vB,t)$, the PDF of $\{\vB_t\}_0^T$ is ,
\eq{\label{kolmo}
\p_t \rho  +\n\cdot(\rho(\vr-\underline{\bm\kappa}\vB -\vu)\star \vB) - \n\cdot\n\cdot[\rho\frac{\sigma^2}2\vB\otimes\vB] = 0,~\rho(\vB,0)=\rho^0(\vB),
}
for all $t\in(0,T)$ and all $\vB\in \RR:={\R^+}^d$. 
The solution of \eqref{reform} is also the solution of
\eq{&& \label{control3}
\min_{\vu\in\U} J(\vu):=\int_{\RR\times(0,T)}\left[|\vB-\vB^d|^2 -\bm\alpha\cdot \vu(\vB,t) + \beta |\sigma\vB\n \vu(\vB,t)|^2\right]\rho(\vB,t)\d B\d t,
}
subject to (\ref{kolmo}).
The conditions  for having equivalence between the two control problems are detailed in \cite{lebris}. 

\subsection{Computation of  gradients}
Consider the variational form of the Kolmogorov equation: find $\rho$ such that, for all $\hat\rho$,
\[
\int_{\RR}\left(\hat\rho\p_t  \rho -\rho(\vr-\underline{\bm\kappa} \vB -\vu)\star\vB \cdot\n\hat\rho + \frac{\sigma^2}2\n\hat\rho\cdot\n\cdot(\vB\otimes\vB\rho)\right) =0,~~\rho(0)=\rho^0.
\]
Calculus of variations says that a variation $\delta\vu$ yields a  $\delta\rho$ with $\delta\rho(0)=0$ and
\eq{\label{del1}
\int_{\RR}\left(\hat\rho\p_t  \delta\rho -\delta\rho(\vr-\underline{\bm\kappa} \vB -\vu)\star\vB \cdot\n\hat\rho + \frac{\sigma^2}2\n\hat\rho\cdot\n\cdot(\vB\otimes\vB\delta\rho)\right) 
=-\int_{\RR} \rho \vB\star\delta\vu\n\hat\rho.
}
Define the adjoint $\rho^*$ by $\rho^*(T)=0$ and, for all $\hat\rho$,
\eq{\label{adjointprime}&\ds
\int_{\RR}\Big(\p_t\rho^*  \hat\rho & + \hat\rho(\vr-\underline{\bm\kappa} \vB -\vu)\star\vB \cdot\n\rho^* - \frac{\sigma^2}2\n\rho^*\cdot\n\cdot(\vB\otimes\vB\hat\rho)
\cr&&
+\hat\rho(|\vB-\vB^d|^2 - \bm\alpha\cdot \vu + \beta |\sigma\vB\n \vu|^2) \Big)=0.
}
Adding \eqref{del1} with $\hat\rho=\rho^*$ to \eqref{adjointprime}  with $\hat \rho=\delta\rho$  gives
\[ 
\int_{\RR}[\p_t(\rho^*\delta\rho) + (|\vB-\vB^d|^2- \bm\alpha\cdot \vu + \beta |\sigma\vB\n \vu|^2)\delta \rho] =-\int_{\RR} \rho \vB\star\delta\vu\n\rho^*.
\]
As $\rho^*(T)=0$ and $\delta\rho(0)=0$, an integration in time gives
\[
 \int_{\RR\times]0,T[}(|\vB-\vB^d|^2- \bm\alpha\cdot \vu + \beta |\sigma\vB\n \vu|^2)\delta\rho =-\int_{\RR\times]0,T[} \rho \vB\star\delta\vu\n\rho^*.
 \]
 Finally,  by differentiating $J$ in (\ref{control3}), 
 \eq{\label{gradj}
 \delta J  &= &\int_{\RR\times[0,T]}\Big[\left(|\vB-\vB^d|^2 - \bm\alpha\cdot \vu +\beta |\sigma\vB\n \vu|^2)\right)\delta \rho 
 \cr
 ~~~~~ &-&\bm\alpha\cdot\rho\delta \vu +2\rho
 \beta \sigma^2\vB\n \vu:\vB\n \delta\vu)\Big]
 \cr&
 =&-\int_{\RR\times]0,T[} \rho \Big[\vB\star\delta\vu\n\rho^*
 +\bm\alpha\cdot\delta \vu -2
 \beta \sigma^2\vB\n \vu:\vB\n \delta\vu\Big]
}
The computation of the gradient follows, because $\delta J=<\hbox{grad}_u J,\delta \vu> +o(|\delta\vu|)$.

\subsection{Numerical Simulation}
Two species are considered (d=2) with $\sigma=\sigma'=0.1$, $\vB_1(0)=1.2$, $\vB_2(0)=0.8$,
\[
\vr=\left[\begin{matrix}1.5 \cr 1.5\end{matrix}\right], 
\underline{\bm\kappa}=\left[\begin{matrix}1.2& -0.1 \cr 0.1 & 1.2\end{matrix}\right],
\bm\alpha=\left[\begin{matrix}0.1 \cr 0.1\end{matrix}\right],
\bm\beta=\left[\begin{matrix}0.02 \cr 0.02\end{matrix}\right], q=1.3, u_m=0.4, u_M=1.4.
\]

A numerical simulation has been done using \texttt{freefem} \cite{freefem}, the finite element method and the optimization module \texttt{ipopt} (see {\small\texttt{https://github.com/coin-or/Ipopt}}). 
Before optimization $J=-0.24$ and after optimization $J=-0.32$.

For simplicity it is assumed that $\vu$ depends on $\vB$ but not on $t$; it was shown numerically in \cite{MLGPOP2} that the time dependence is small.

 The main difficulty is due to the non integrability of the right hand side in the adjoint equation. At all levels $\RR$ must be replaced by a finite domain smaller than the infinite integration domain of the partial differential equations. 
Results are shown on the following 4 figures. Figure \ref{resDP2}, \ref{kolmofig} show the surfaces
$\vu_i$, i=1,2, functions of $\vB_1,\vB_2$.

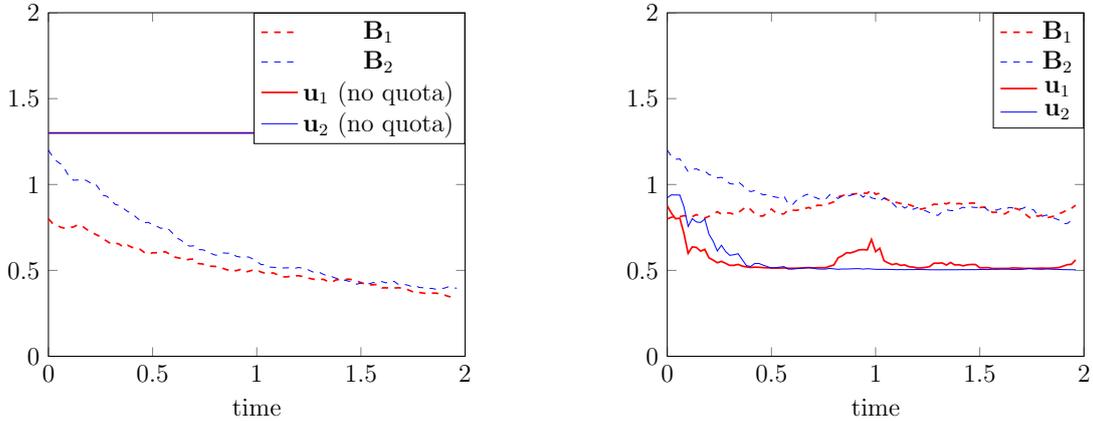
\begin{figure}[h!]
\begin{minipage} [b]{0.49\textwidth}
\begin{center}
\begin{tikzpicture}[scale=0.8]
\begin{axis}[legend style={at={(1,1)},anchor=north east}, compat=1.3,
  xmin=0, xmax=2,ymin=0,ymax=2,
  xlabel= {time}]
\addplot[thick,dashed,color=red,mark=none,mark size=1pt] table [x index=0, y index=1]{fig/nocontrol.txt};
\addlegendentry{ $\vB_1$}
\addplot[thin,dashed,color=blue,mark=none,mark size=1pt] table [x index=0, y index=2]{fig/nocontrol.txt};
\addlegendentry{ $\vB_2$}
\addplot[thick,solid,color=red,mark=none,mark size=1pt] table [x index=0, y index=3]{fig/nocontrol.txt};
\addlegendentry{ $\vu_1$ (no quota)}
\addplot[thin,solid,color=blue,mark=none,mark size=1pt] table [x index=0, y index=4]{fig/nocontrol.txt};
\addlegendentry{ $\vu_2$ (no quota)}
\end{axis}
\end{tikzpicture}
\end{center}
\end{minipage}
\hskip 0.25cm
\begin{minipage} [b]{0.49\textwidth}
\begin{center}
\begin{tikzpicture}[scale=0.8]
\begin{axis}[legend style={at={(1,1)},anchor=north east}, compat=1.3,
  xmin=0, xmax=2,ymin=0,ymax=2,
  xlabel= {time}]
\addplot[thick,dashed,color=red,mark=none,mark size=1pt] table [x index=0, y index=1]{fig/optcontrol.txt};
\addlegendentry{ $\vB_1$}
\addplot[thin,dashed,color=blue,mark=none,mark size=1pt] table [x index=0, y index=2]{fig/optcontrol.txt};
\addlegendentry{ $\vB_2$}
\addplot[thick,solid,color=red,mark=none,mark size=1pt] table [x index=0, y index=3]{fig/optcontrol.txt};
\addlegendentry{ $\vu_1$}
\addplot[thin,solid,color=blue,mark=none,mark size=1pt] table [x index=0, y index=4]{fig/optcontrol.txt};
\addlegendentry{ $\vu_2$}
\end{axis}
\end{tikzpicture}
\end{center}
\end{minipage}
\caption{Sample trajectories computed with the  control from the Kolmogorov equation with and without quota for 2 species. The corresponding quotas are also shown. Without quota the biomass decays with time dangerously.
\label{price} }
\end{figure}

With this optimal quota, two sample trajectories where chosen randomly.  Results are shown on Figure \ref{price}.  Similar trajectories without quota are given for comparison on the left.

\begin{figure}[h!]
\begin{minipage} [b]{0.49\textwidth}
\centering
\begin{tikzpicture}[scale=0.9]
\begin{axis}[ legend style={at={(1,1)},anchor=north east}, compat=1.3,xlabel= {$\vB_1$},ylabel= {$\vB_2$}]
 \addplot3[surf,fill opacity=0.75] table [ ] {fig/u1.txt};
\addlegendentry{ $\vu_1(\vB_1,\vB_2)$}
\end{axis}
\end{tikzpicture}
\caption{Kolmogorov solution: value of $\vu_1(\vB)$.}
\label{resDP2} 
\end{minipage}
\begin{minipage} [b]{0.49\textwidth}
\centering
\begin{tikzpicture}[scale=0.9]
\begin{axis}[legend style={at={(1,1)},anchor=north east}, compat=1.3,xlabel= {$\vB_1$},ylabel= {$\vB_2$}]
 \addplot3[surf,fill opacity=0.5] table [ ] {fig/u2.txt};
\addlegendentry{  $\vu_2(\vB_1,\vB_2)$}
\end{axis}
\end{tikzpicture}
\caption{Kolmogorov solution: value of $\vu_2(\vB)$.}
\label{kolmofig}
\end{minipage}
\end{figure}

\section{Quotas Computed by a Markovian Neural Network}

Here too, let us simplify the problem by forgetting the time dependency  of the quota and represent each component of $\vu(\vB)$ by a Neural Network with $K=2$ hidden layers of 50 neurons each and \texttt{ReLU} activations. Denote $\vX=(\vB,\vu)^T$, so that the NN represents also $\vB\mapsto \vu(\vB)$:
\[
\vX^0 \hbox{ given },  \vX^{k+1} := \sum _{k=0}^K\max\{\underline{\vA}^k \vX^k + \vb^k,0\},~~ \vu_{NN}(\vB):=\left[\begin{matrix}I & 0\end{matrix}\right]\vX^K.
\]
Then the coefficients $\underline{\vA^k}$ and $\vb^k$ are computed by minimizing $J$ (the `loss') defined by \eqref{reform} with $\vu_{NN}$ in place of $\vu$.

This method was proposed and tested in \cite{MLGPOP2} but Pierre Bras \cite{PBGP} gave a convergence proof when a modified version (called Langevin) of the stochastic optimization algorithm ADAM is used.  For the numerical tests we used his open source implementation with \texttt{Keras} (see {\small \texttt{https://github.com/Bras-P/langevin-for-stochastic-control}}).

\begin{figure}[h!]
\begin{minipage} [b]{0.29\textwidth}
\begin{center}
\begin{tikzpicture}[scale=0.65]
\begin{axis}[legend style={at={(1,1)},anchor=north east}, compat=1.3,
  xmin=-0.05,
  xlabel= {epoch}]
\addplot[thick,solid,color=red,mark=none,mark size=1pt] table [x index=1, y index=2]{fig/loss.txt};
\addlegendentry{loss static}
\addplot[thick,solid,color=blue,mark=none,mark size=1pt] table [x index=1, y index=2]{fig/lossT.txt};
\addlegendentry{loss dynamic}
\addplot[thin,dashed,color=red,mark=none,mark size=1pt] table [x index=1, y index=4]{fig/loss.txt};
\addplot[thick,dashed,color=blue,mark=none,mark size=1pt] table [x index=1, y index=5]{fig/loss.txt};
\addplot[thin,dashed,color=red,mark=none,mark size=1pt] table [x index=1, y index=4]{fig/lossT.txt};
\addplot[thick,dashed,color=blue,mark=none,mark size=1pt] table [x index=1, y index=5]{fig/lossT.txt};
\end{axis}
\end{tikzpicture}
\end{center}
\end{minipage}
\hskip 0.25cm
\begin{minipage} [b]{0.29\textwidth}
\begin{center}
\begin{tikzpicture}[scale=0.65]
\begin{axis}[legend style={at={(1,1)},anchor=north east}, compat=1.3,
  xmin=0, xmax=2,ymin=0,ymax=2,
  xlabel= {time}]
\addplot[thick,dashed,color=red,mark=none,mark size=1pt] table [x index=0, y index=1]{fig/traj.txt};
\addlegendentry{ $\vB_1$}
\addplot[thick,dashed,color=blue,mark=none,mark size=1pt] table [x index=0, y index=2]{fig/traj.txt};
\addlegendentry{ $\vB_2$}
\addplot[thin,solid,color=red,mark=none,mark size=1pt] table [x index=0, y index=3]{fig/traj.txt};
\addlegendentry{ $\vu_1$ static}
\addplot[thin,solid,color=blue,mark=none,mark size=1pt] table [x index=0, y index=4]{fig/traj.txt};
\addlegendentry{ $\vu_2$ static}
\end{axis}
\end{tikzpicture}
\end{center}
\end{minipage}
\hskip 0.25cm
\begin{minipage} [b]{0.29\textwidth}
\begin{center}
\begin{tikzpicture}[scale=0.65]
\begin{axis}[legend style={at={(1,1)},anchor=north east}, compat=1.3,
  xmin=0, xmax=2,ymin=0,ymax=2,
  xlabel= {time}]
\addplot[thick,dashed,color=red,mark=none,mark size=1pt] table [x index=0, y index=1]{fig/trajT.txt};
\addlegendentry{ $\vB_1$}
\addplot[thick,dashed,color=blue,mark=none,mark size=1pt] table [x index=0, y index=2]{fig/trajT.txt};
\addlegendentry{ $\vB_2$}
\addplot[thin,solid,color=red,mark=none,mark size=1pt] table [x index=0, y index=3]{fig/trajT.txt};
\addlegendentry{ $\vu_1$ dynamic}
\addplot[thin,solid,color=blue,mark=none,mark size=1pt] table [x index=0, y index=4]{fig/trajT.txt};
\addlegendentry{ $\vu_2$ dynamic}
\end{axis}
\end{tikzpicture}
\end{center}
\end{minipage}
\caption{Solution of the  control problem with 2 species using $\vu(\vB)$ (called static) compared with using $\vu(\vB,t)$ (dynamic). The loss functions are displayed on the left. After optimization of the loss (left) $J=-0.0315$ in the static case and 0.15 in the dynamic case.  In the middled (static) and on the right (dynamic) two sample trajectories (blue and orange) and their control (green and red) are displayed. The results with dynamic controls are poor.}
\label{bras2NN}
\end{figure}
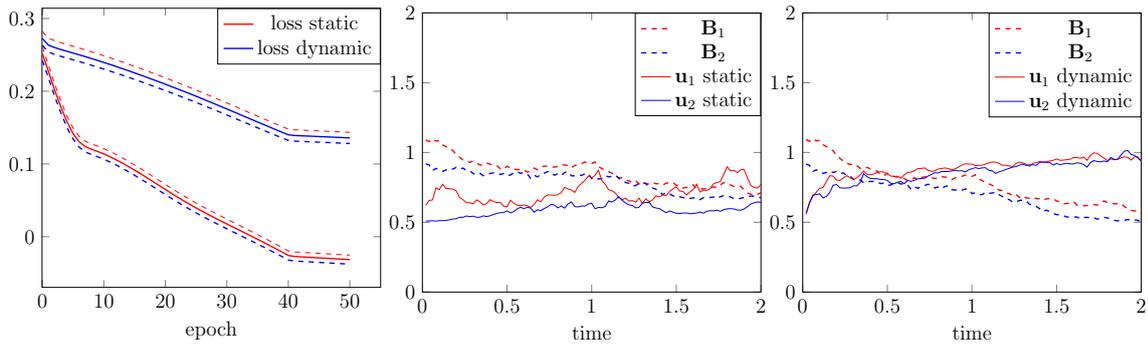

The numerical results are shown on Figure \ref{bras2NN} on the same problem described above.
The converged value of the loss function is greater than the Kolmogorov solution which is typical because Supervised Learning does not compute the absolute minimum but on the other hand the solution proposed is usually more robust.

The biggest asset of Supervised Learning is that it can be used with any number of species while Dynamic Programming cannot be used beyond 3 species.

\section{A Simple Strategy}
\begin{figure}[h!]
\begin{minipage}[b]{0.5\textwidth}
\begin{center}
\begin{tikzpicture}[scale=0.8]
\begin{axis}[legend style={at={(1,1)},anchor=north east}, compat=1.3,
  xmin=0, xmax=2,ymin=0,ymax=2,
  xlabel= {time}]
\addplot[thick,solid,color=red,mark=none,mark size=1pt] table [x index=0, y index=1]{fig/easyu2D.txt};
\addlegendentry{ $\vu_1$}
\addplot[thin,dashed,color=red,mark=none,mark size=1pt] table [x index=0, y index=2]{fig/easyu2D.txt};
\addlegendentry{ $\vB_1$}
\addplot[thick,solid,color=blue,mark=none,mark size=1pt] table [x index=0, y index=3]{fig/easyu2D.txt};
\addlegendentry{ $\vu_2$}
\addplot[thin,dashed,color=blue,mark=none,mark size=1pt] table [x index=0, y index=4]{fig/easyu2D.txt};
\addlegendentry{ $\vB_2$}
\end{axis}
\end{tikzpicture}
\end{center}
\caption{Stabilization of the biomass of 2 species by the simple control of \eqref{easyc}. The constraints $\vu_1,\vu_2\in[0.4,1.4]$ do not break the method in this case.}
\label{easyu}
\end{minipage}
\hskip 0.25cm
\begin{minipage} [b]{0.5\textwidth}
\begin{center}
\includegraphics[width=6cm]{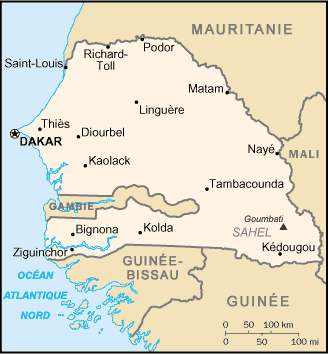}
\caption{Sketched map of Senegal (Wikipedia).}
\end{center}
\label{senegal}
\end{minipage}
\end{figure}

Common sense tells us that if the biomass is decreasing (resp. increasing) then the quota should be made smaller (resp. bigger). In practical terms this means
\eq{\label{easyc}
\vu(t+\delta t)=\vu(t) +\omega(\vB(t)-\vB(t-\delta t)).
}
Figure \ref{easyu} shows the results for the same problem as above with $\omega=100$. This simple solution may stabilize the biomasses at their initial levels but it cannot bring them to a desire level different from the initial value. Furthermore, it does not account for the political cost of the quota, $\bm\alpha\cdot\vu$.


\section{A Fishing Model with Quotas in the Open Sea}

\subsection{A Behavioral model for fishermen}
All variables are now  function of spatial $\vx$ and time $t$. 
Recognizing that $\nabla\vB$ is a local indicator for a better fishing site, the position of a fishing boat $\vZ(t)$ is driven by
\eq{
\dot \vZ(t) =U_M \frac{\nabla \vB}{|\nabla \vB|}|_{\vZ(t),t},~~B(0)=B^0.
}
where $U_M$ is the cruise speed of the boat.
To be profitable the amount of fish caught should be greater that the  operating cost, itself proportional to the square of the velocity of the boat, i.e.
 \eq{\label{profit}
 \bm\gamma\cdot{\vB}(\vZ(t),t) >U_M^2, \hbox{ otherwise the fisherman returns home}.
 }

\subsection{The Logistic equation for the Biomass}
Assume that fish move with a  velocity $\vv$ and a small randomness $\nu$.   The velocity $\vv$ could be the sea current  plus their own velocity to follow the plankton gradient $\nabla P$ where $P$ is the plankton biomass.  

Fishing  depletes the fish population as before but only where fishing occurs. So if $M$ is the number of boats, then  at point $\vx$ of the domain studied $\Omega\subset\R^2$, and time $t$, the  fish biomass $\vB(\vx,t)$ is driven by a PDE in $\Omega\times(0,T)$,
\eq{\label{fish}
\p_t{\vB} + \nabla\cdot(\vv(\nabla P) \vB)- \nu\Delta{\vB} = \vB \star\left(P\vr -\sum_1^M\vu(\vZ^i,t) - \underline{\bm\kappa} \vB \right), \quad \vB(0)={\vB}^0,
}
with  $\partial \vB/\partial\vn=0$ on the border $\Gamma$ of $\Omega$ where $\vn$ is its outer normal to $\Omega$.  Plankton contributes to the reproductive welfare of fish by a positive factor for each species $P\vr $. The total catch is $\vB\star\vu$; as before $\kappa$ is the capacity matrix of the site. In practice it is strongly dependent on $\vx$ but in absence of  information we ran the model with $\kappa$ constant. 

The problem is 
\eq{\label{opt2x}
\ds \min_{\vu\in \U} J:=\int_0^T\left(\int_\Omega|\vB(t)-\vB^d(t)|^2-\sum_1^M(\bm\alpha\cdot\vu(\vZ^i,t)-\bm\beta[\vu](\vZ^i,t))\right)\d t
}
subject to \eqref{fish}

\begin{remark}
It may be feasible to replace \eqref{fish} by a system  equivalent at the limit $\delta t\to 0$:
\eq{\label{fishstoc}&
\vB(\vx,t)&=\vB(\vx-\vv(\vx,t)\delta t,t-\delta t) \cr&&
 + \delta t \vB\star\left[P\vr-\sum_1^M\vu(\vZ^i,t) - \underline{\bm\kappa} \vB \right]_{|\vx,t-\delta t} + 2\sqrt{\nu\delta t}\vN_0^1
, \hbox{ for all $\vx$,}
 }
   The long time limit could be studied with the stationary Kolmogorov equation for the invariant measure of the process. 
\end{remark}

\subsection{A Logistic Equation for the Plankton}
Letting the fish drift with the currents is too simple. If fish follows a plankton density ${P}$ then $\vv\nabla \vB$ in \eqref{fish} is replaced by $\n\cdot(\vB \n P)$.
Assume plankton is regenerated at rate one and eaten by some fish species at rate $\vb $. The logistic equation for $P$ is:
\eq{\label{plankton}
\p_t{P} + \vv\cdot\n{P}-\mu\Delta {P} =  {P}(1- P - \vb\cdot{ \vB}),~~ \frac{\p P}{\p n}|_\Gamma=0~or ~ P|_\Gamma=0,~~{P}(0)={P}^0
}
where $P^0(x)$ is the plankton density  at initial time. The model assumes that in absence of fish the long time limit (the fishing site plankton capacity)  of $P$ is one.
Here $\vv$ is the sea current velocity. Other models, perhaps more realistic, can be found in \cite{franks}.
 
\begin{remark} 
{\rm If $\vb\cdot \vB <1$, then ${P}$ is positive and bounded by $1-\vb\cdot \vB$, if it is initially so.  Otherwise ${P}$ may become negative and the model is no longer meaningful.
}\end{remark}
\begin{remark}
{\rm When $c':=1-\vb\cdot \vB$ is constant and $\vv=0$ and $\mu=0$, the solution of $\dot{P}={P}(c'-{P})$ is ${P}=c' \e^{c' t} /(1+\e^{c't})$, and it tends to $c'$ when $t\to\infty$.
When $\mu>0,~\vv=0$ and $\Omega$ is bounded, then   limit $\lim_{t\to \infty}P=c'$}.
\end{remark}

\subsection{Numerical Simulation Without Quota}

We ran the model with one species only but with plankton, with $\Omega$ a portion of the Atlantic Ocean facing Senegal (see Figure \ref{senegal}), with the following parameters,
\[
 T=2, \delta t=0.02, c=0.7, a=0.2,b=1, \mu=0.1,r=1, \kappa=1, 
 K=100, U_M=2,  \gamma=1. 
\]
 A random noise of variance $ \sigma=0.05$ is added to the position of the boats at each time iteration. Initialization is
\[
 Q_{t=0}=0.05,~~P^0 = [1-\frac1{40}((x-4)^2+(y-6)^2)]^+ ,~ 
B^0 = [1-\frac1{40}((x-4)^2+(y-6)^2)]^+.
\]

To obtain a meaningful sea current we set
\[
		\vv = 10\cos(2\pi t)\n\psi \hbox{ where }
\Delta\psi =0,~~\psi|_{\Gamma_1}=\vx_1-6,~~\psi|_{\Gamma_2}=0.
\]
where $\Gamma_1$ and $\Gamma_2$ are the upper and lower boundaries of the domain.

\begin{figure}[htbp]
\begin{center}
\includegraphics[width=7.8cm]{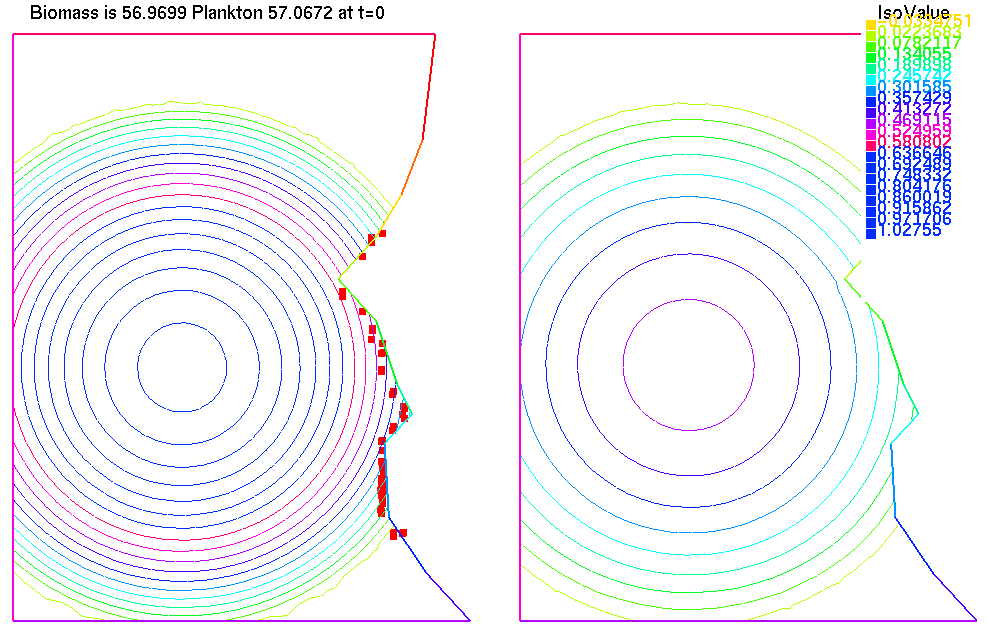}
\includegraphics[width=7.8cm]{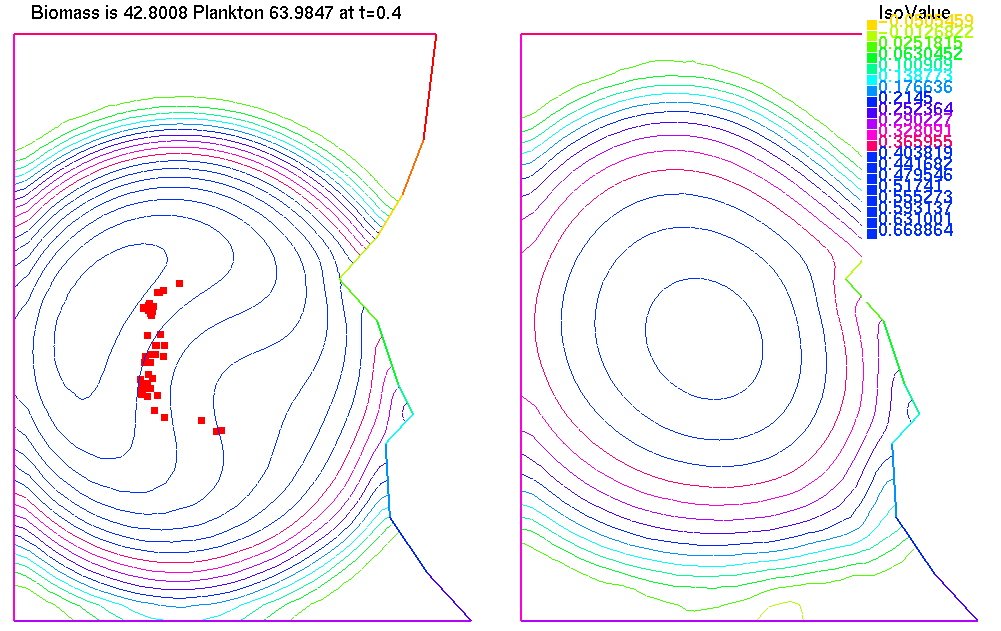}

\includegraphics[width=7.8cm]{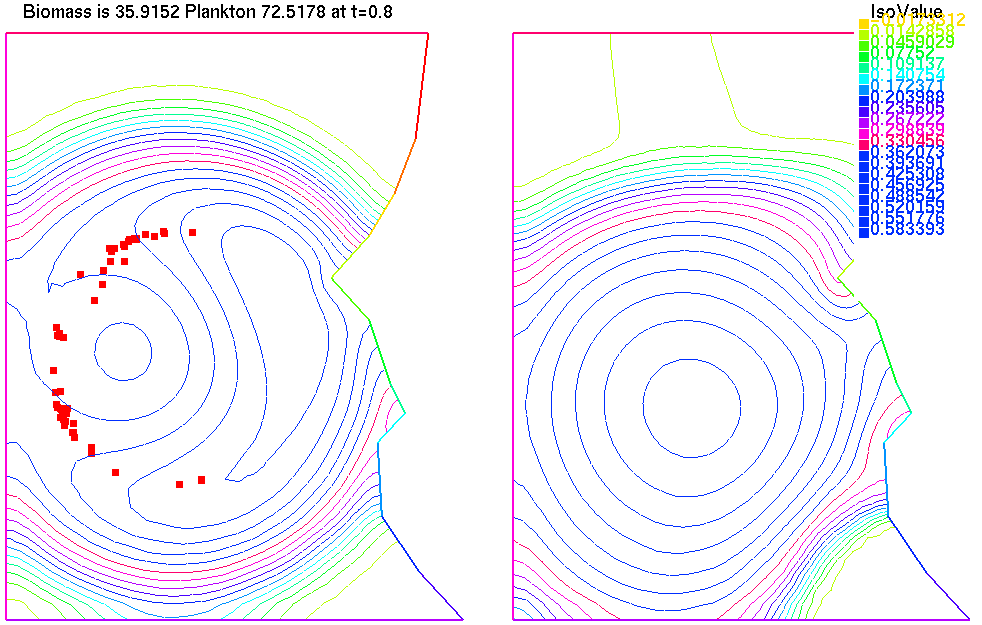}
\includegraphics[width=7.8cm]{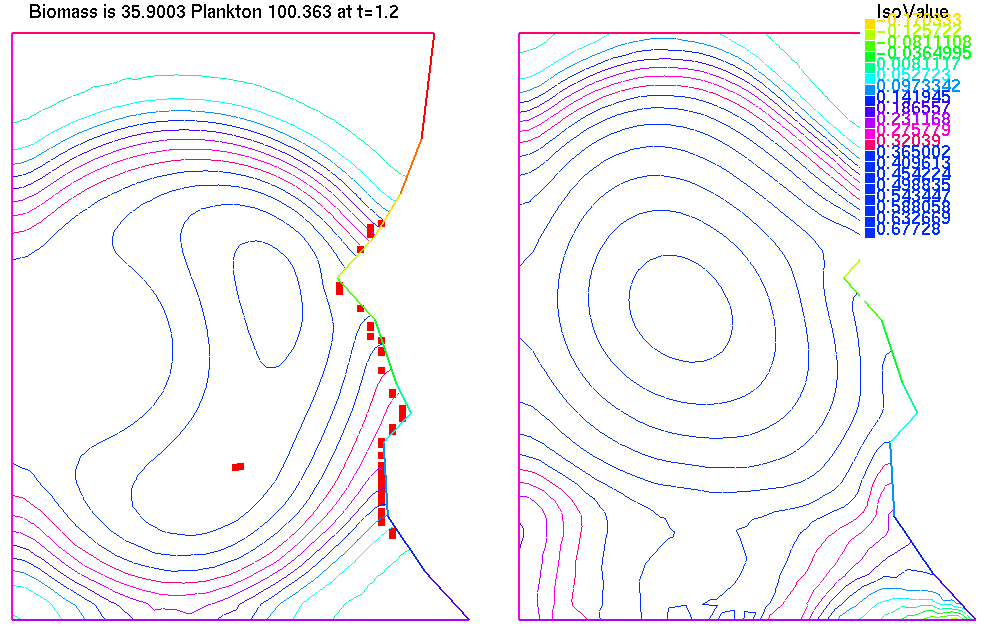}

\caption{  From left to right and top to bottom:  Level lines of fish (left) and plankton (right) biomass at times $t=0.,0.4,0.8, 1.2$. The color map legends apply to $B$ and $P/2$. The total biomass and plankton are indicated above the B-plots. The positions of the 50 fishing boats are indicated by small red squares. In this case \underline{without quota} the fishermen fish extensively until $t=1$ and then run out of resource (fishing is no longer profitable) and go back to the coast.  }
\label{simnoquota}
\end{center}
\end{figure}

The following plots in Figure \ref{simnoquota} show 1/ the initial position of the 50 boats on the coast and the level lines of $B$ (left) and $2P$ (right), 2/ their position and the values of $B$ and $P$ at time at 0.4, then  3/, 4/ are the same but at time 0.8 and 1.2.  The integrals of $P$ and $B$ in $\Omega$ are displayed on top of the plots of $B$ and also on Figure \ref{bilan}..

We see  that the fishing boats move towards the maximum zone of $B$ and then spread because the biomass reduces drastically. Shortly after $t=1$ the catch is too small for profit (see \eqref{profit}) so the boats return to the coast and stay there until $t=T$.

\subsection{Numerical Simulation with Quota}
 All parameters are as above but now $Q$ is adjusted by 
 \eq{\label{quota}
 Q(t+\delta t)=Q(t)+\delta t\int_\Omega (B_{t}-B_{t-\delta t})\d x.
}
We see on Figure \ref{simwithquota} that the behavior is very different with quota. The boats move to the maximum zone of $B$ but stay there because the quota prevents to fishermen from depleting the biomass.  The boats stay at the same spot till $B$ plateaus and  the boat positions spread due to the noise added to $\vZ$ at each time step.

This is seen too on Figure  \ref{bilan} which shows the evolution with time of the mean of $B$, the mean of $P$ and the mean of $Q$.  
\begin{figure}[htbp]
\begin{center}
\includegraphics[width=7.8cm]{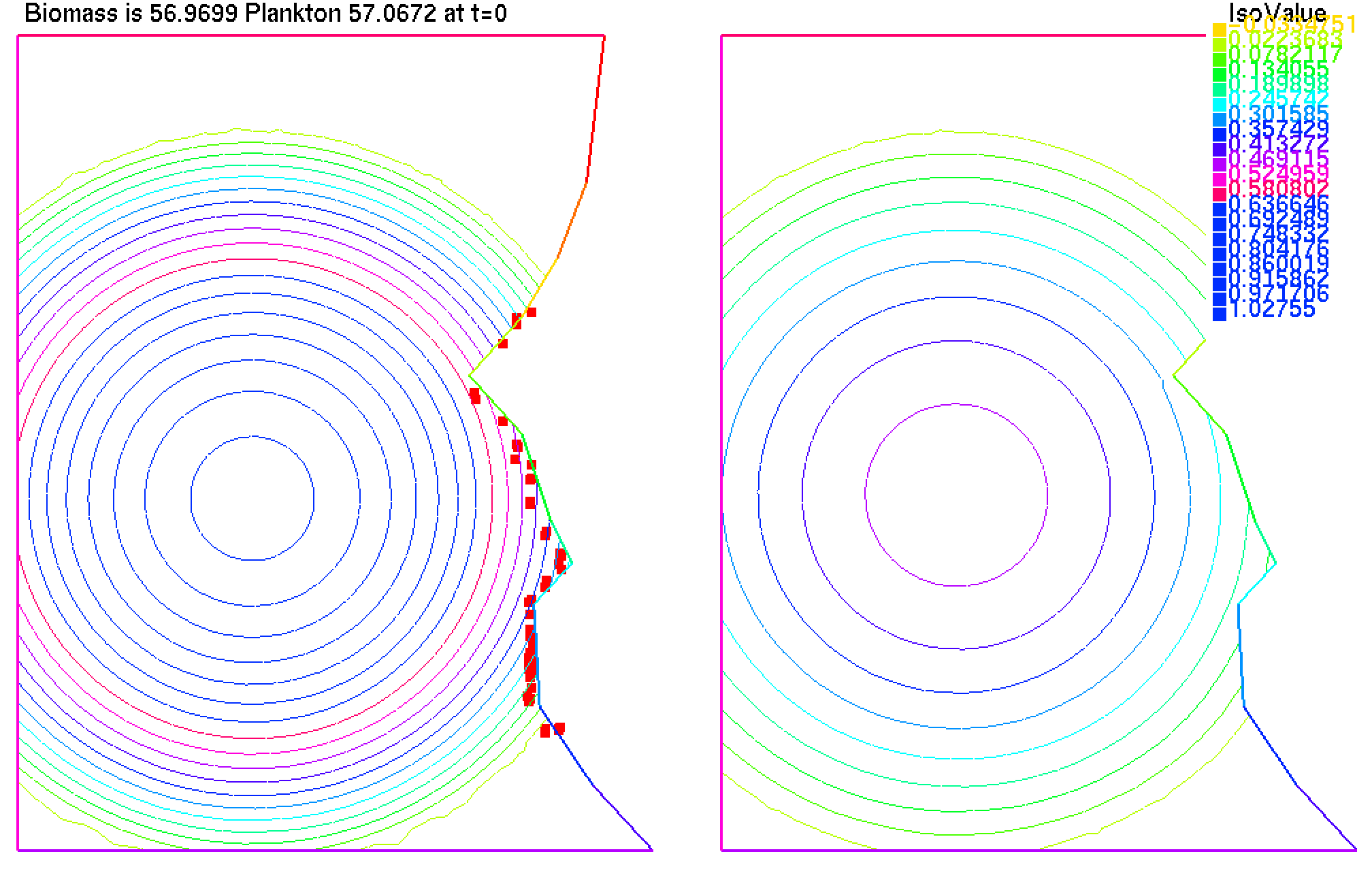}
\includegraphics[width=7.8cm]{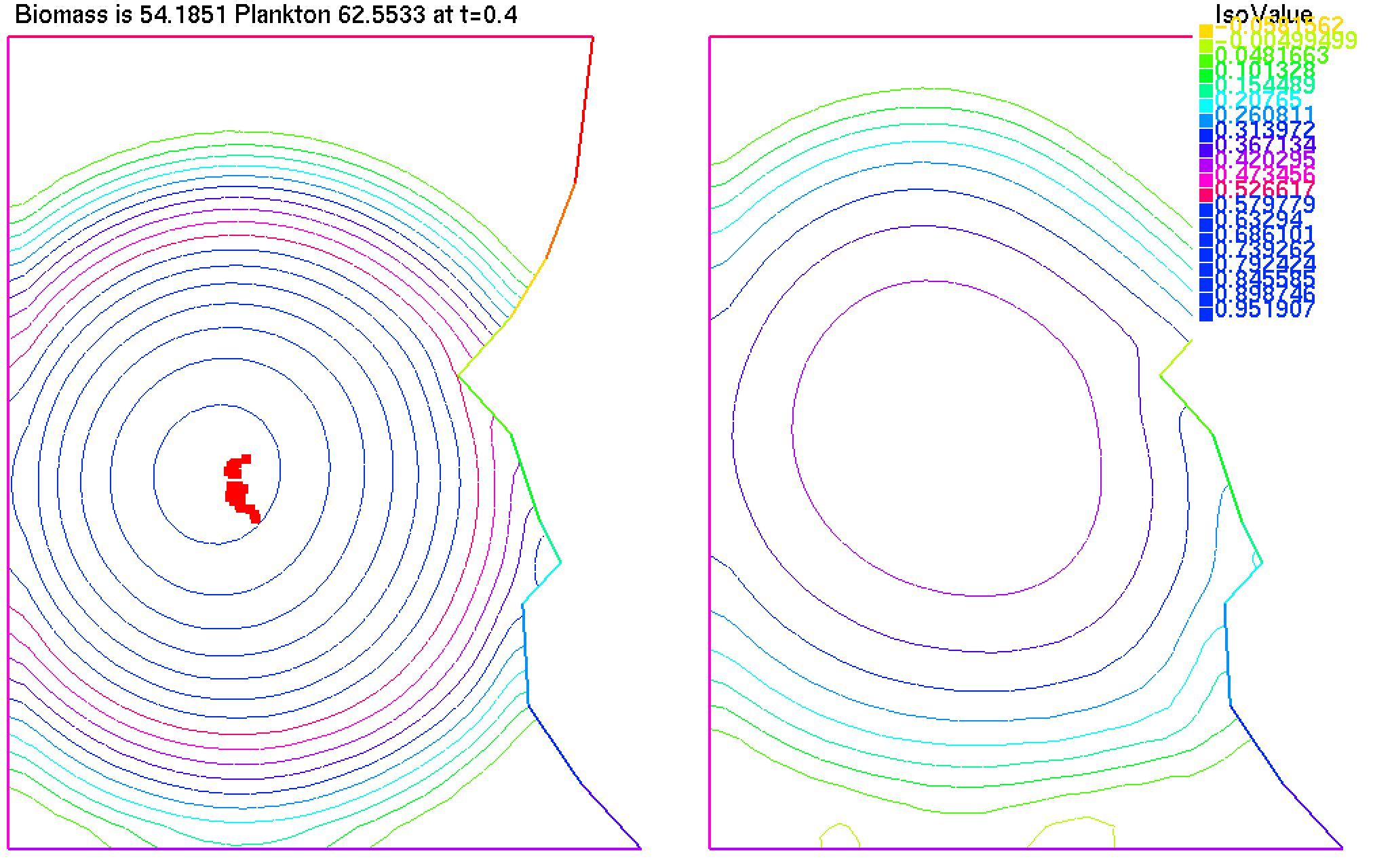}
\includegraphics[width=7.8cm]{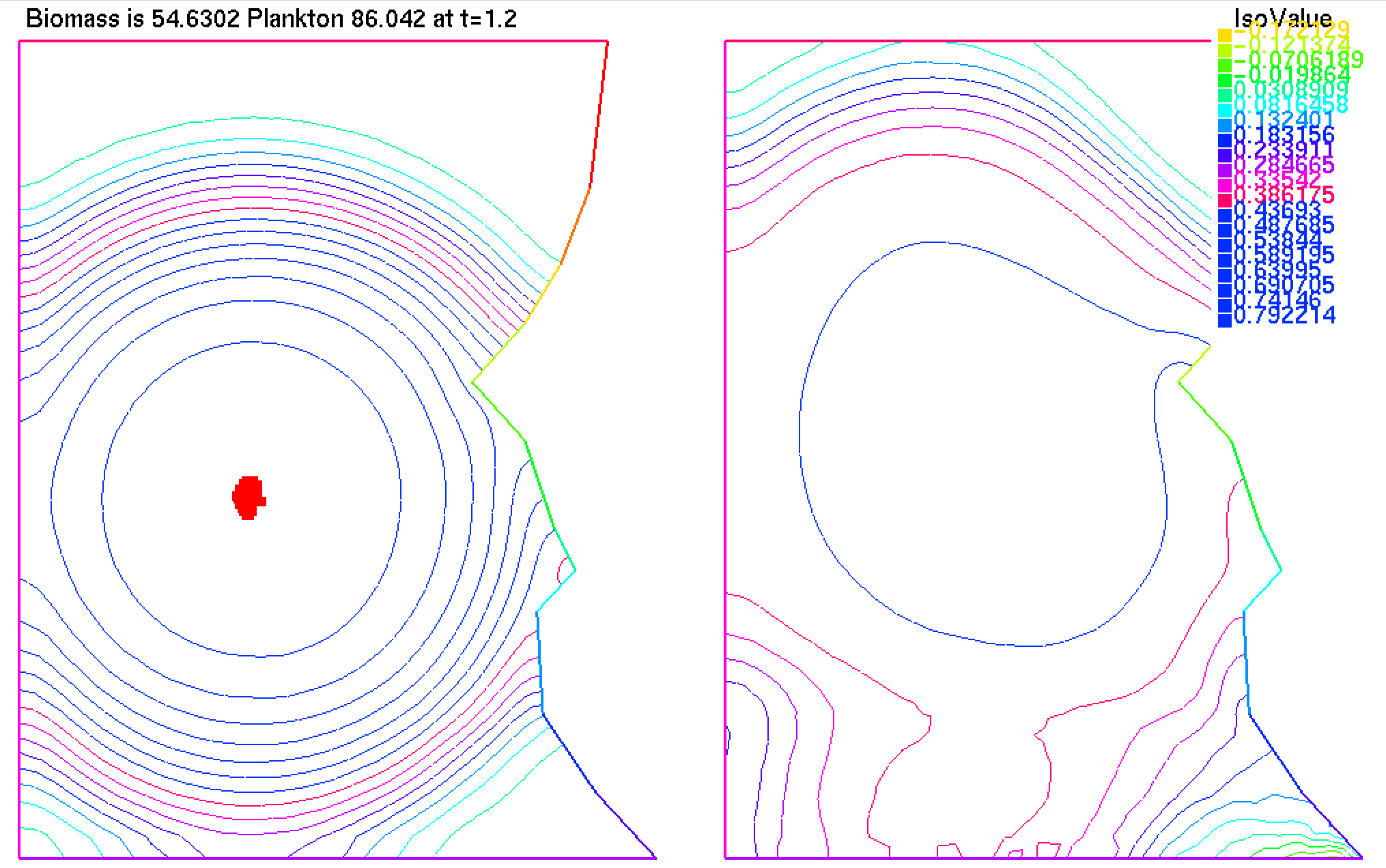}
\includegraphics[width=7.8cm]{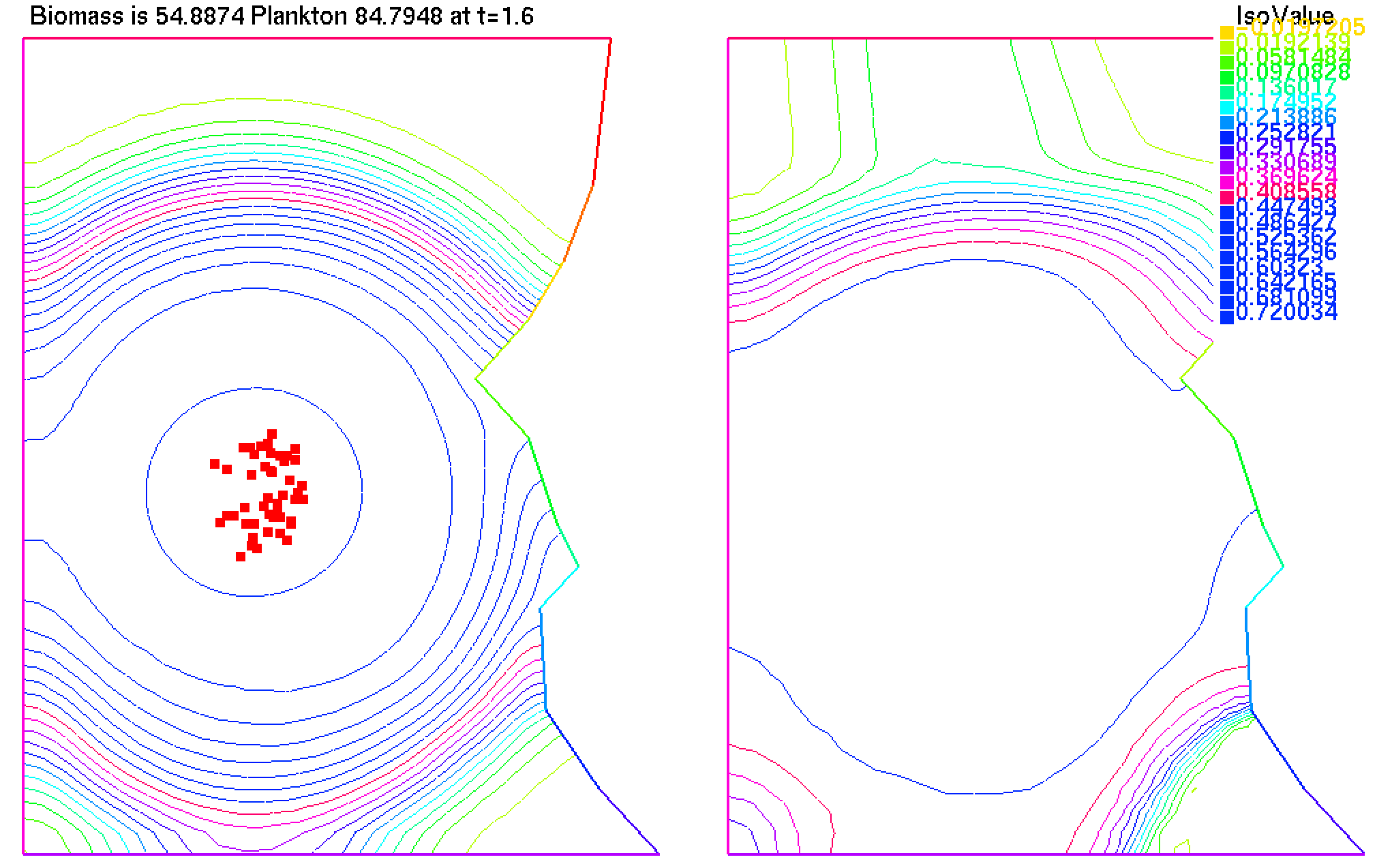}
\caption{  From left to right and top to bottom:  Level lines of fish (left) and plankton (right) biomass at  times $0.,0.4,1.2, 1.6$. The color map legends apply to $B$ and $P/2$. The total biomass and plankton are indicated above the B-plots. The positions of the 50 fishing boats are indicated by small red squares. In this case \underline{with quota} the fishermen sail to the maximum of the biomass but as the catch is limited by the quota, $B$ stays above the level of profitability at all time. Later $B$ plateaus over a large area in the center of the domain and so the fishermen to not correct the spatial scattering due to the noise.  }
\label{simwithquota}
\end{center}
\end{figure}
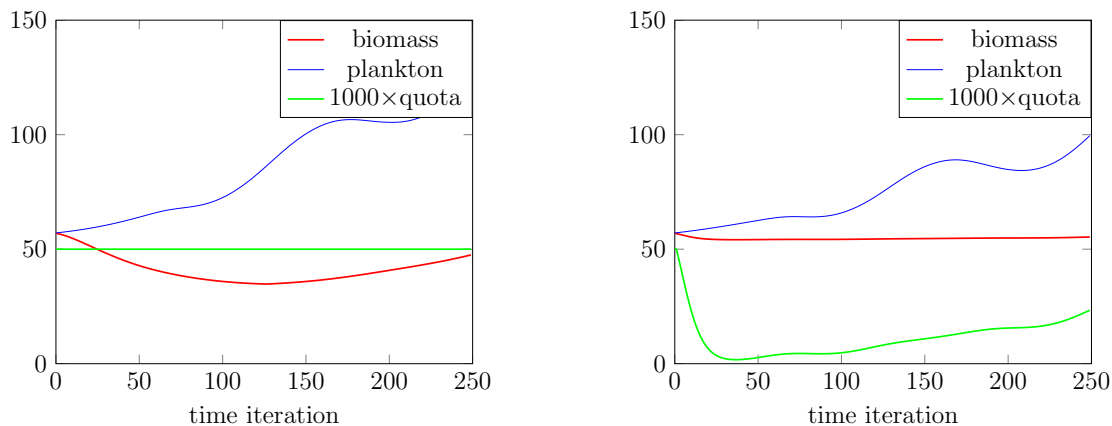
\begin{figure}[h!]
\begin{minipage} [b]{0.49\textwidth}
\begin{center}
\begin{tikzpicture}[scale=0.8]
\begin{axis}[legend style={at={(1,1)},anchor=north east}, compat=1.3,
  xmin=0, xmax=250,ymin=0,ymax=150,
  xlabel= {time iteration}]
\addplot[thick,solid,color=red,mark=none,mark size=1pt] table [x index=0, y index=1]{fig/histnq.txt};
\addlegendentry{ biomass}
\addplot[thin,solid,color=blue,mark=none,mark size=1pt] table [x index=0, y index=2]{fig/histnq.txt};
\addlegendentry{ plankton}
\addplot[thick,solid,color=green,mark=none,mark size=1pt] table [x index=0, y index=3]{fig/histnq.txt};
\addlegendentry{ $1000\times $quota}
%
\end{axis}
\end{tikzpicture}
\end{center}
\end{minipage}
\hskip 0.25cm
\begin{minipage} [b]{0.49\textwidth}
\begin{center}
\begin{tikzpicture}[scale=0.8]
\begin{axis}[legend style={at={(1,1)},anchor=north east}, compat=1.3,
  xmin=0, xmax=250,ymin=0,ymax=150,
  xlabel= {time iteration}]
\addplot[thick,solid,color=red,mark=none,mark size=1pt] table [x index=0, y index=1]{fig/histwq.txt};
\addlegendentry{ biomass}
\addplot[thin,solid,color=blue,mark=none,mark size=1pt] table [x index=0, y index=2]{fig/histwq.txt};
\addlegendentry{ plankton}
\addplot[thick,solid,color=green,mark=none,mark size=1pt] table [x index=0, y index=3]{fig/histwq.txt};
\addlegendentry{ $1000\times $quota}
%
\end{axis}
\end{tikzpicture}
\end{center}
\end{minipage}
\caption{Evolution of the total biomass $\int_\Omega B$ and scaled total plankton $\int_\Omega P/2$ with and without quota. Notice (on the right) that the quota strategy \eqref{quota} is very efficient at maintaining the biomass constant. The quotas are displayed in green, it is constant by hypothesis on the left figure. }
\label{bilan}
\end{figure}

\section*{Conclusion}
With the single site model of \cite{MLGPOP2}, we have confronted two methods to adjust the quotas for single sites fisheries and shown that Supervised Learning does fairly well on a problem with 2 species.  For more than 2 species only Supervise Learning is applicable. Then we have put some foundation stones for a distributed model for fishing in the Atlantic ocean facing Senegal and shown that a common sense strategy to keep the biomass constant works. We have seen that the effect of quotas on the fishing strategy of fishermen is striking. A more sophisticated strategy is yet to be found for the control of the biomasses in large areas like the Atlantic ocean.  Whatever has been said for fisheries translates to several other population control problems but once again these are theoretical case studies which are far from applicable directly to real life situations. 

\section*{Acknowledgement}
We thank P. Auger and M. Lauri\`ere for their helpful comments;  All PDE computations have been done with the public domain FreeFEM++ \cite{freefem}.

\bibliographystyle{plain}
\bibliography{biblio}

\begin{thebibliography}{10}

\bibitem{ALL}
P.~M. Allen and J.~M. McGlade.
\newblock Modelling complex human systems: A fisheries example.
\newblock {\em European Journal of Operational Research}, 30:147--167, 1987.

\bibitem{PAOP}
P.~Auger and O.~Pironneau.
\newblock {\em Parameter Identification by Statistical Learning of a Stochastic
  Dynamical System Modelling a Fishery with price variation}.
\newblock Comptes rendus de l'acad\'emie des sciences, 2020.

\bibitem{BIC}
A.~Bick.
\newblock Quadratic-variation-based dynamic strategies.
\newblock {\em Management Sciences}, 41(4):722--732, 1995.

\bibitem{PBGP}
P.~Bras and G.~Pag\`es.
\newblock Convergence of langevin-simulated annealing algorithms with
  multiplicative noise ii: Total variation.
\newblock {\em Monte Carlo Methods and Applications},
  doi:10.1515/mcma-2023-2009, 2023.

\bibitem{lebris}
C.~Le Bris and P.~L. Lions.
\newblock Existence and uniqueness of solutions to fokker-planck type equations
  with irregular coefficients.
\newblock {\em Comm}, 33:1272--1317, 2008.

\bibitem{MP18}
T.~Brochier, P.~Auger, D.~Thiao, A.~Bah, S.~Ly, T.~Nguyen Huu, and P.~Brehmer.
\newblock Can overexploited fisheries recover by self-organization?
  reallocation of the fishing effort as an emergent form of governance.
\newblock {\em Marine Biology}, 95:46--56, Mar 2018.

\bibitem{DAN}
A.~Danielsson.
\newblock Efficiency of catch and effort quotas in the presence of risk.
\newblock {\em Journal of Environmental Economics and Management}, 43:20--33,
  2002.

\bibitem{mcglade}
J.~McGlade (ed).
\newblock The dynamics of flows of matter and energy.
\newblock {\em L. Pimm, Chapter 6, Blackwell Science}, 6, 1999.

\bibitem{haddon}
M.~Haddon.
\newblock Modelling and quantitative methods in fisheries, crc press, taylor \&
  francis.
\newblock 2011.

\bibitem{freefem}
F.~Hecht.
\newblock New development in freefem++, j.
\newblock {\em Numer. Math.}, 20:251--265, 2012.

\bibitem{katsukawa}
T.~Katsukawa.
\newblock Numerical investigation of the optimal control rule for
  decision-making in fisheries management.
\newblock {\em Fisheries Science}, 70:123--131, 2004.

\bibitem{king}
M.~King.
\newblock Ecology and ecosystem in fisheries biology, assessment and
  management, blackwell publishing.
\newblock 1995.

\bibitem{MLGPOP2}
M.~Lauri\`ere, G.~Pag\`es, and O.~Pironneau.
\newblock Performance of a markovian neural network versus dynamic programming
  on a fishing control problem.
\newblock {\em Probability, Uncertainty and Quantitative Risk}, 8(1):121--140,
  2023.

\bibitem{moussaoui}
A.~Moussaoui, M.~Bensenane, P.~Auger, and A.~Bah.
\newblock On the optimal size and number of reserves in a multi-site fishery
  model.
\newblock {\em Journal of Biological Systems}, 23(01):31--47, 2015.

\bibitem{franks}
Peter J. S.~Franks: Npz.
\newblock models of planton dynamics.
\newblock {\em J. of Oceanography}, 58:379--387, 2002.

\bibitem{pagespironneau}
G.~Pag\`es and O.~Pironneau.
\newblock Protection of a fishing site with optimal quotas.
\newblock to appear, 2020.

\bibitem{Punt}
A.~Punt, D.~Butterworth, C.~deMoor, J.~DeOliveira, and M.~Haddon.
\newblock Management strategy evaluation: best practices.
\newblock {\em Fish and Fisheries}, John Wiley(DOI: 10.1111/faf.12104), 2014.

\bibitem{schaefer}
M.~B. Schaefer.
\newblock Some aspects of the dynamics of populations important to the
  management of commercial marine fisheries.
\newblock {\em Inter-American Tropical Tuna Commission}, 2526, 1954.

\bibitem{verhulst}
P.-F. Verhulst.
\newblock Notice sur la loi que la population poursuit dans son accroissement.
\newblock {\em Correspondance math\'ematique et physique, No}, 10:113--121,
  1838.

\bibitem{lotka}
V.~Volterra.
\newblock Variations and fluctuations of the number of individuals in animal
  species living together.
\newblock In R.~N Chapman, editor, {\em Animal}. Ecology. McGraw--Hill., 1931.

\bibitem{WEI}
M.~Weitzman.
\newblock Landing fees vs harvest quotas with uncertain fish stocks.
\newblock {\em Journal of Environmental Economics and Management}, 43:325--338,
  2002.

\end{thebibliography}

\end{document}